\documentclass[a4paper,11pt]{article}

\usepackage[top=1in, bottom=1in, left=1in, right=1in]{geometry}

\usepackage[utf8]{inputenc}
\usepackage[T1]{fontenc}
\usepackage[british]{babel}
\usepackage{csquotes}

\usepackage{lmodern}

\usepackage{amsthm}

\theoremstyle{plain}
\newtheorem{theorem}{Theorem}
\newtheorem{lemma}[theorem]{Lemma}
\newtheorem{corollary}[theorem]{Corollary}

\theoremstyle{definition}
\newtheorem{definition}[theorem]{Definition}

\theoremstyle{remark}
\newtheorem{remark}[theorem]{Remark}

\usepackage[backend=biber,maxbibnames=99]{biblatex}
\usepackage{type1cm}        %
\usepackage{makeidx}         %
\usepackage{graphicx}        %
\usepackage{multicol}        %
\usepackage[bottom]{footmisc}%

\usepackage{newtxtext}       %
\usepackage[varvw]{newtxmath}       %

\usepackage{tikz}
\usetikzlibrary{patterns}

\usepackage{dsfont}

\usepackage{mathtools}

\usepackage{subcaption}

\usepackage{booktabs}

\newcommand{\bbR}{\mathbb{R}}

\newcommand{\calH}{\mathcal{H}}
\newcommand{\calT}{\mathcal{T}}

\DeclareMathOperator{\supp}{supp}

\newcommand{\matr}[1]{\mathbf{#1}}
\renewcommand{\vec}[1]{\mathbf{#1}}

\title{A GenEO-type coarse space with smaller eigenproblems}
\author{Peter Bastian and Nils Friess$^\ast$}
\date{{\footnotesize $^\ast$Interdisciplinary Center for Scientific Computing, Heidelberg University, Germany\\
 \texttt{\{peter.bastian,nils.friess\}@iwr.uni-heidelberg.de}}%
 }

\begin{document}
\maketitle

\abstract{Coarse spaces are essential to ensure robustness w.r.t.\ the number of subdomains in two-level overlapping Schwarz methods. Robustness with respect to the coefficients of the underlying partial differential equation (PDE) can be achieved by adaptive (or spectral) coarse spaces involving the solution of local eigenproblems. The solution of these eigenproblems, although scalable, entails a large setup cost which may exceed the cost for the iteration phase. In this paper we present and analyse a new variant of the GenEO (\underline{Gen}eralised \underline{E}igenproblems in the \underline{O}verlap) coarse space which involves solving eigenproblems only in a strip connected to the boundary of the subdomain. This leads to a significant reduction of the setup cost while the method satisfies a similar coefficient-robust condition number estimate as the original method, albeit with a possibly larger coarse space.}

\section{Introduction and problem setting}
\label{sec:1}
Coarse spaces are essential to ensure robustness w.r.t.\ the number of subdomains in two-level overlapping Schwarz methods. Robustness with respect to the coefficients of the underlying partial differential equation (PDE) can be achieved by adaptive (or spectral) coarse spaces involving local eigenproblems \cite{spillane_geneo_2014,heinlein_adaptive_2019,bastian_multilevel_2023}. The solution of these eigenproblems, although scalable, entails a large setup cost which may exceed the cost for the iteration phase. Following an idea from \cite{msgfem_ring} we present and analyse a new variant of the GenEO (\underline{Gen}eralised \underline{E}igenproblems in the \underline{O}verlap) coarse space~\cite{spillane_geneo_2014} which solves eigenproblems only in a strip connected to the boundary of the subdomain. If the overlap parameter is not too large the setup cost is significantly reduced. The method satisfies a coefficient-robust condition number estimate similar to that of the original method, but the size of the coarse space might be increased. We call the new method R-GenEO as the domain of the eigenproblems has ring-shape in two dimensions. Another interpretation might be ``real'' GenEO as the method now really only solves eigenproblems in the overlap, in contrast to the original method.

In the rest of this section we introduce the problem setting and review the original GenEO coarse space. In Section~\ref{sec:rgeneo} we introduce the new R-GenEO coarse space and prove its robustness.
In the last section we give a brief numerical comparison of both methods.

\subsection{Preliminaries}
We consider a finite element discretisation of a partial differential equation that is posed on some domain $\Omega \subset \bbR^d$, $d=2,3$. Let $V_h$ be a finite-dimensional Hilbert space (the finite element space on some mesh $\calT_h$), $a:V_h \times V_h \rightarrow \bbR$ a continuous, symmetric and coercive bilinear form and $f \in V_h^\prime$ a linear form. We assume that $a$ is of the form 
\begin{equation*}
    a(u,v) = \sum_{\tau \in \calT_h} a_\tau(u_{|\tau}, v_{|\tau})\,, \qquad \text{ for all } u,v \in V_h\,,
\end{equation*}
for some positive semidefinite bilinear forms $a_\tau : V_h(\tau) \times V_h(\tau) \rightarrow \bbR$, for all $\tau \in \calT_h$, where $V_h(\tau) = \{ v_{|\tau} \mid v \in V_h\}$. The problem then reads: find $u \in V_h$ such that $a(u,v) = \langle f, v \rangle$ for all $v \in V_h$. While the method is in principle more general, we restrict ourselves to a scalar, second-order elliptic problem where $a(u,v) = \int_\Omega (\alpha(x) \nabla u) \cdot \nabla v \, \mathrm{d}x$. By choosing a basis $\{\varphi_k\}_{k=1,\dotsc,n}$ of $V_h$ we can write this problem as a linear system $\matr{A}\vec{u}=\vec{f}$. Our goal is to solve this linear system using a Krylov method preconditioned by a two-level Schwarz method. 

To this end, let $\{ \Omega_j \}_{j=1}^N$ be an overlapping domain decomposition of the domain $\Omega$ where we assume that the subdomains are resolved by the mesh $\calT_h$, and introduce the spaces $V_{h,0}(\Omega_j) \coloneqq \{ v_{|\Omega_j} \mid v \in V_h\,, \supp(v) \subset \Omega_j \}$. We denote by $R_j^\top : V_{h,0}(\Omega_j) \rightarrow V_h$ the extension-by-zero operator. Its adjoint $R_j : V_h^\prime \rightarrow V_{h,0}(\Omega_j)^\prime$ is called the restriction operator. In addition we introduce a so-called coarse space $V_H \subset V_h$ which will be defined later. We denote by $R_H^\top : V_H \rightarrow V_h$ the natural embedding and by $R_H$ its adjoint. If we let $\matr{R}_j$, $j=1,\dotsc,N$, be the matrix representations of $R_j$ (w.r.t.\ the basis $\{ \varphi_k \}_k$) and $\matr{R}_H$ be the matrix representation of $R_H$ (w.r.t.\ a basis of the coarse space $V_H$) then the additive two-level Schwarz preconditioner reads
\begin{equation*}
    \matr{M}_{\mathrm{AS},2}^{-1} = \matr{R}_H^T \matr{A}_H^{-1} \matr{R}_H + \sum_{j=1}^N \matr{R}_j^T \matr{A}_j^{-1} \matr{R}_j\,,
\end{equation*}
where the subdomain matrices $\matr{A}_j$ and the coarse matrix $\matr{A}_H$ are defined via Galerkin projection, i.e., $\matr{A}_j = \matr{R}_j \matr{A} \matr{R}_j^T$ and  $\matr{A}_H = \matr{R}_H \matr{A} \matr{R}_H^T$.

Let us recall some standard definitions and results. We denote by $k_0 \in \mathbb{N}$ the maximum number of subdomains that a grid cell belongs to, i.e.,  
\begin{equation*}
    k_0 = \max_{\tau \in \calT_h} \left( \#\left\{ \Omega_j \mid j \in \{1, \dotsc, N\}, \tau \in \Omega_j \right\} \right)\,.
\end{equation*}

For any set $D$ that is the union of elements of $\calT_h$ we let $V_h(D) \coloneqq \{ v_{|D} \mid v \in V_h \}$ and write
\begin{equation*}
    a_D(u,v) \coloneqq \sum_{\tau \in D} a_\tau(u_{|\tau}, v_{|\tau}) \qquad \text{ for } u,v \in V_h(D)\,.
\end{equation*}

Note that for $u, v \in V_{h,0}(\Omega_j)$ we have $a_{\Omega_j}(u, v) = a(R_j^\top u, R_j^\top v)$ and this bilinear form is positive definite (on $V_{h,0}(\Omega_j)$). We denote by $|\cdot|_{a,D}$ the seminorm induced by the bilinear form $a_{D}(\cdot, \cdot)$, and by $\|\cdot\|_{a,\Omega_j}$ the norm induced by $a_{\Omega_j}(\cdot, \cdot)$ on $V_{h,0}(\Omega_j)$. If $D = \Omega$, we omit the domain in the subscript.

\begin{definition}[Stable decomposition]
    \label{def:stable decomposition}
    Let $C_0 > 0$ be a constant. A $C_0$-stable decomposition of $u \in V_h$ is a family of functions $\{z_j\}_{j=0,\dotsc,N}$ such that 
    \begin{equation}
        \label{eq:stable decomposition 1}
        u = \sum_{j=0}^N z_j\,, \qquad \text{with } z_0 \in V_H \text{ and } z_j \in V_{h,0}(\Omega_j) \text{, for } j > 0\,,
    \end{equation}
    and
    \begin{equation}
        \label{eq:stable decomposition 2}
        \| z_0 \|_a^2 + \sum_{j=1}^N \|z_j\|^2_{a, \Omega_j} \le C_0^2 \| u \|_a^2\,.
    \end{equation}
\end{definition}

Proving existence of a stable decomposition is the key step in deriving condition number bounds for additive Schwarz methods. The following result appears in some form or another in many publications on additive Schwarz methods.

\begin{theorem}[{{\cite[Thm.\ 2.8]{spillane_geneo_2014}}}]
    If every $u \in V_h$ admits a $C_0$-stable decomposition then 
    $\kappa(\matr{M}_{\mathrm{AS},2}^{-1} \matr{A}) \le C_0^2(k_0 + 1)$.
\end{theorem}

At the expense of introducing a quadratic dependency on $k_0$ in the condition number bound, one can show~\eqref{eq:stable decomposition 2} without having to bound $\|z_0\|_a$ in terms of $\|u\|_a$. 

\begin{lemma}[{{\cite[Lem.\ 2.9]{spillane_geneo_2014}}}]
    Using the notation of Definition~\ref{def:stable decomposition}, if there exists $C_1 > 0$ s.t.
    \begin{equation*}
        \| z_j \|_{a, \Omega_j}^2 \le C_1 | u |^2_{a, \Omega_j}\,,\qquad \text{ for all } j = 1, \dotsc, N\,,
    \end{equation*}
    then the decomposition~\eqref{eq:stable decomposition 1} is $C_0$-stable with $C_0^2 = 2 + C_1 k_0(2k_0 + 1)$.
\end{lemma}

\subsection{The GenEO coarse space}
Let us now briefly recall the definition of the original GenEO coarse space introduced in~\cite{spillane_geneo_2014}. To this end, we first define the overlapping zone (or just the overlap, for short) of subdomain $\Omega_j$, $j = 1, \dotsc, N$, as
\begin{equation*}
    \Omega^\circ_j = \left\{ x \in \Omega_j \mid x \in \Omega_{j^\prime} \text{ for some } j^\prime \neq j \right\}\,.
\end{equation*}
Next, let $\xi_j : V_h(\Omega_j) \rightarrow V_{h,0}(\Omega_j)$ be partition of unity operators that satisfy
\begin{equation*}
    \sum_{j=1}^N R_j^\top \xi_j(v_{|\Omega_j}) = v\,, \qquad \text{ for all } v \in V_h\,,
\end{equation*}
and $\xi_j(v)_{|\Omega_j \setminus \Omega_j^\circ} = v_{|\Omega_j \setminus \Omega_j^\circ}$, for all $v \in V_h(\Omega_j)$ and $j = 1, \dotsc, N$. On each overlapping subdomain $\Omega_j$, $j=1, \dotsc, N$, we now consider the generalised eigenvalue problem: Find $(\lambda^j, t^j) \in \bbR \times V_h(\Omega_j)$ such that
\begin{equation}
  \label{eq:eigenproblem geneo}
  a_{\Omega_j}(u, t^j) = \lambda^j a_{\Omega_j^\circ}(\xi_j(u), \xi_j(t^j))\,,\qquad \text{for all $u \in V_h(\Omega_j)$}.
\end{equation}

For $j = 1, \dotsc, N$ let $\{ t_j^k \}_{k = 1, \dotsc, m_j}$ be the eigenfunctions of the eigenproblem~\eqref{eq:eigenproblem geneo} corresponding to the $m_j$ smallest eigenvalues. The GenEO coarse space is defined as
\begin{equation*}
    V_H \coloneqq \mathrm{span} \left\{
        R^\top_j \xi_j( t^j_k ) \mid k=1,\dotsc,m_j\,; j = 1, \dotsc, N
    \right\}\,.
\end{equation*}

In~\cite{spillane_geneo_2014} the authors prove that the condition number of the lwo-level Schwarz method with the GenEO coarse space can be bounded by
\begin{equation*}
    \kappa(\matr{M}^{-1}_{\mathrm{AS}, 2} \matr{A}) \le (1 + k_0) ( 2 + k_0(2 k_0 + 1) \max_{1 \le j \le N} (1 + 1 / \lambda^j_{m_j + 1} ) )\,.
\end{equation*}

\section{The R-GenEO coarse space}\label{sec:rgeneo}
To define the new coarse space we introduce some additional sets and notation. Let $\Omega_j^\ast \subseteq \Omega_j$ be a set that satisfies $\Omega_j^\circ \subset \Omega_j^\ast$ (e.g.\ $\Omega_j^\ast$ can be obtained by extending $\Omega_j^\circ$ by one layer of grid elements towards the interior of $\Omega_j$). We further define
\begin{align*}
    \omega_j^\circ \coloneqq \overline{\Omega_j} \setminus \overline{\Omega_j^\circ}\,,
    &&
    \Gamma_j^\circ \coloneqq \partial \Omega_j^\circ \,\cap\, \partial \omega_j^\circ\,,
    &&
    \omega_j^\ast \coloneqq \overline{\Omega_j} \setminus \overline{\Omega_j^\ast}\,,
    &&
    \Gamma_j^\ast \coloneqq \partial \Omega_j^\ast \,\cap\, \partial \omega_j^\ast\,,
\end{align*}
and we assume that there exists $\delta > 0$ such that $\text{dist}(\Gamma_j^\ast, \Gamma_j^\circ) \ge \delta$.
Now consider the following eigenvalue problem on each overlapping subdomain: Find $(\lambda^j, t^j) \in \bbR \times V_h(\Omega_j^\ast)$ such that
\begin{equation}
  \label{eq:eigenproblem}
  a_{\Omega_j^\ast}(u, t^j) = \lambda^j a_{\Omega_j^\ast}(\eta_j(u), \eta_j(t^j))\,,\qquad \text{for all $u \in V_h(\Omega_j^\ast)$},
\end{equation}
where $\eta_j : V_h(\Omega_j^\ast) \rightarrow V_{h,0}(\Omega_j^\ast)$ are functions that satisfy $(\eta_j v)(x) = (\xi_j v)(x)$ for $x \in \overline{\Omega_j^\circ}$, in particular $(\eta_jv)(x)=0$ for $x\in\Gamma_j^\ast$ and $(\eta_jv)(x)=1$ for $x \in \Gamma_j^\circ$.
Let $t_k^j$, $k = 1, \dotsc, m_j$, denote the eigenvectors corresponding to the $m_j$ smallest eigenvalues of~\eqref{eq:eigenproblem}. To obtain the local components of the basis vectors $y_k^j$ that make up the coarse space we proceed as follows: on $\overline{\Omega_j^\circ}$ we take $y_k^j = t_k^j$. On $\omega_j^\circ$ we compute $y_k^j$ as an operator-harmonic extension of $t_k^j|_{\Gamma_j^\circ}$. More precisely, consider the $a$-orthogonal projection operator $Q_j^\circ : V_h(\omega_j^\circ) \rightarrow V_{h,0}(\omega_j^\circ)$ which is defined via
\begin{equation*}
    a_{\omega_j^\circ}(Q_j^\circ v, z) = a_{\omega_j^\circ}(v, z) \qquad \text{for all $z \in V_{h,0}(\omega_j^\circ)$}\,,
\end{equation*}
and set $H_j^\circ = I - Q_j^\circ$. To define the harmonic extension of the eigenfunctions (which are elements of $V_h(\Omega_j^\ast)$) we introduce the trace operator $\gamma_j^\circ$ defined by $\gamma_j^\circ v = v|_{\Gamma_j^\circ}$ and an extension operator $E_j^\circ : V_h(\Gamma_j^\circ) \rightarrow V_h(\omega_j^\circ)$  that satisfies $\gamma_j^\circ E_j^\circ = I$ on $V_h(\Gamma_j^\circ)$ (i.e., $ E_j^\circ$ takes a finite element function on the boundary $\Gamma_j^\circ$ and extends it arbitrarily to the interior domain $\omega_j^\circ$). The local components of the coarse space vectors are then defined as
\begin{equation}
  \label{eq:coarse vectors}
  y_k^j(x) =
  \begin{cases}
    t_k^j(x) \,, &\text{for  $x \in \overline{\Omega_j^\circ}$,} \\
    H_j^\circ E_j^\circ \gamma_j^\circ t^j_k(x) \,, &\text{for $x \in \omega_j^\circ$,}
  \end{cases}
\end{equation}
for $k = 1, \dotsc, m_j$. We introduce the short notation $y_k^j = \calH_{\Omega_j^\circ \rightarrow \Omega_j}(t_k^j)$ for the mapping defined by~\eqref{eq:coarse vectors}. In practice, the second line of~\eqref{eq:coarse vectors} amounts to discarding the values of $t_k^j$ in $\Omega_j^\ast \cap \omega_j^\circ$ and computing the harmonic extension from $\Gamma_j^\circ$ to $\omega_j^\circ$.

\begin{definition}[R-GenEO coarse space]
  For $j = 1, \dotsc, N$ let $\{ t_j^k \}_{k = 1, \dotsc, m_j}$ be the eigenfunctions of the eigenproblem~\eqref{eq:eigenproblem} corresponding to the $m_j$ smallest eigenvalues. The R-GenEO coarse space is defined as
  \begin{equation*}
    V_H \coloneqq \text{span} \left\{
      R^\top_j \xi_j( \calH_{\Omega_j^\circ \rightarrow \Omega_j}(t^j_k) ) \mid k=1,\dotsc,m_j\,; j = 1, \dotsc, N
    \right\}\,.
  \end{equation*}
\end{definition}

\begin{remark}
   As mentioned above, a similar construction was recently used in the context of a multiscale generalised finite element method~\cite{msgfem_ring}. The idea of computing the coarse basis vectors using energy-minimising extensions of eigenvectors is also present in the so called adaptive GDSW coarse spaces, see, e.g.,~\cite{heinlein_adaptive_2019}. 
\end{remark}

To prove robustness of the coarse space, we make use of the following result which was shown in~\cite{bastian_multilevel_2023}.

\begin{lemma}[{{\cite[Sec.\ 3.3]{bastian_multilevel_2023}}}]
    \label{lem:stability general}
    Let $V$ be a $n$-dimensional vector space and let $a,b : V \times V \rightarrow \bbR$ be two positive semidefinite bilinear forms on $V$ with $\ker a \cap \ker b = \{ 0 \}$. Consider the generalised eigenvalue problem: Find $(\lambda, p) \in (\bbR \cup \{\infty\}) \times V$, $p \neq 0$, such that either $p \not\in \ker b$ and 
    \begin{equation*}
        a(p, v) = \lambda b(p, v) \qquad \text{ for all } v \in V
    \end{equation*}
    or $p \in \ker b$ and $\lambda = \infty$. Let the eigenpairs $\{ (p_k, \lambda_k)\}_{k=1}^n$ of this problem be ordered such that $0 \le \lambda_1 \le \dotsm \le \lambda_n \le \infty$. Suppose that $m \in \{1, \dotsc, n\}$ is such that $0 < \lambda_{m+1} < \infty$. Then, the projection operator
    \begin{equation*}
        \Pi_{m} v \coloneqq \sum_{k = 1}^{m} b(v, p_k) p_k
    \end{equation*}
    is well-defined and orthogonal w.r.t.\ the bilinear form $a(\cdot, \cdot)$. Thus ${| \Pi_m v |}_{a} \le {|v|}_{a}$ and ${| v - \Pi_m v |}_{a} \le {|v|}_{a}$ and we have the stability estimate
    \begin{equation*}
      {| v - \Pi_{m} v |}^2_{b} \le 1 / \lambda_{m + 1} {| v - \Pi_{m} v |}_{a}^2
    \end{equation*}
    for all $v \in V$.
\end{lemma}

To apply the lemma, we have to show that the kernels of the bilinear forms that appear in the eigenproblem~\eqref{eq:eigenproblem} have trivial intersection. The proof is the same as for the classical GenEO coarse space (see~\cite[Lem.\ 3.18]{bastian_multilevel_2023}).

\begin{lemma}
    For $j \in \{1, \dotsc, N\}$ let $a(\cdot, \cdot) = a_{\Omega_j^\ast}(\cdot, \cdot)$ and $b(\cdot, \cdot) = a_{\Omega_j^\ast}(\eta_j(\cdot), \eta_j(\cdot))$. Then $\ker a \cap \ker b = \{ 0 \}$.
\end{lemma}
\begin{proof}
    For subdomains where the extended overlapping zone $\Omega_j^\ast$ touches the global Dirichlet boundary, $a(\cdot, \cdot)$ is positive definite, hence $\ker a = \{ 0 \}$. Otherwise, $\ker a = \mathrm{span} \{ \mathds{1} \} $, where $\mathds 1$ is the constant one function on $\Omega_j$. But $\mathds 1 \not \in \ker b$ due to the modified partition of unity function $\eta_j$.
\end{proof}

Thus we can apply Lemma~\ref{lem:stability general} to the eigenproblem~\eqref{eq:eigenproblem}.

\begin{corollary}
    Let $m_j \in \{1, \dotsc, \dim(V_h(\Omega^\ast_j))\}$ be such that $0 < \lambda^j_{m_j+1} < \infty$. The local projection operator $\Pi^\ast_{j, m_j} : V_h(\Omega_j^\ast) \rightarrow V_h(\Omega_j^\ast)$,
    \begin{equation*}
      \Pi^\ast_{j, m_j} v \coloneqq \sum_{k = 1}^{m_j} a_{\Omega_j^\ast}(\eta_j(v), \eta_j(t^j_k)) t^j_k\,,
    \end{equation*}
    is well-defined and orthogonal w.r.t.\ the bilinear form $a_{\Omega_j^\ast}(\cdot, \cdot)$. Thus
    \begin{equation*}
    {| \Pi^\ast_{j,m_j} v |}_{a, \Omega_j^\ast} \le {|v|}_{a, \Omega_j^\ast} \quad\text{and}\quad{| v - \Pi^\ast_{j,m_j} v |}_{a, \Omega_j^\ast} \le {|v|}_{a, \Omega_j^\ast}    
    \end{equation*}
    and we have the local stability estimate
    \begin{equation*}
      | \eta_j(v - \Pi^\ast_{j,m_j} v) |^2_{a, \Omega_j^\ast} \le 1 / \lambda^j_{m_j + 1} | v - \Pi^\ast_{j,m_j} v |_{a, \Omega_j^\ast}^2\,.
    \end{equation*}
\end{corollary}

The operator $\Pi^\ast_{j, m_j}$ only maps to $V_h(\Omega_j^\ast)$ and as such only defines the coarse components in $\Omega_j^\ast$. We thus additionally define
\begin{equation*}
    \Pi_{j, m_j} v \coloneqq \sum_{k = 1}^{m_j} a_{\Omega_j^\ast}(\eta_j(v), \eta_j(t^j_k)) y^j_k\,.
\end{equation*}

We can now define a stable decomposition.
\begin{theorem}
  Let $v \in V_h(\Omega)$. The decomposition
  \begin{equation*}
    z_0 \coloneqq \sum_{j = 1}^N \xi_j (\Pi_{j,m_j} v|_{\Omega_j})\,,
    \qquad
    z_j \coloneqq \xi_j (v|_{\Omega_j} - \Pi_{j,m_j} v|_{\Omega_j}), \quad \text{for $j = 1, \dotsc, N$}\,,
  \end{equation*}
  is $C_0$-stable with
  \begin{equation*}
    C_0^2 = 2 + k_0(2 k_0 + 1) \max_{1 \le j \le N} (2 + 3 / \lambda^j_{m_j + 1})\,.
  \end{equation*}
\end{theorem}
\begin{proof}
  Since $\xi_j$ is the identity for restrictions of functions to $\omega_j^\circ$ we have
  \begin{equation}
    \label{eq:splitting}
    {|z_j|}_{a, \Omega_j}^2 = {\left|\xi_j (v - \Pi_{j,m_j} v)\right|}_{a, \Omega_j^\circ}^2 + {\left| v - \Pi_{j,m_j} v \right|}_{a, \omega_j^\circ}^2\,.
  \end{equation}
  We will treat each term separately. For the first part, we use that $\xi_j = \eta_j$ and $t^j_k = y^j_k$ on $\Omega_j^\circ$ so that by Lemma~\ref{lem:stability general} we have
  \begin{align*}
      {\left|\xi_j (v - \Pi_{j,m_j} v)\right|}_{a, \Omega_j^\circ}^2 
      &= {\big|\eta_j (v - \Pi^\ast_{j,m_j} v)\big|}_{a, \Omega_j^\circ}^2
      \le {\big|\eta_j (v - \Pi^\ast_{j,m_j} v)\big|}_{a, \Omega_j^\ast}^2 \\
      &\le 1 / \lambda^j_{m_j+1} |v|_{a,\Omega_j^\ast}^2 \le 1 / \lambda^j_{m_j+1} |v|_{a,\Omega_j}^2\,.
  \end{align*}
  For the second part, we first write
  \begin{align*}
      {\big| v - \Pi_{j,m_j} v \big|}_{a, \omega_j^\circ}^2 
      &= {\big| (I - H_j^\circ) v + H_j^\circ E_j^\circ \gamma_j^\circ v - H_j^\circ E_j^\circ \gamma_j^\circ \Pi^\ast_{j,m_j} v \big|}_{a, \omega_j^\circ}^2\,.
  \end{align*}
  Here we added and subtracted $H_j^\circ v$ and used the linearity of the mapping $H_j^\circ \circ E_j^\circ \circ \gamma_j^\circ$. Using the triangle inequality we then get
  \begin{align*}
    {\big| v - \Pi_{j,m_j} v \big|}_{a, \omega_j^\circ}^2 
    &\le 2 {\big| (I - H_j^\circ) v\big|}^2_{a,\omega_j^\circ} + 2 {\big| H_j^\circ E_j^\circ \gamma_j^\circ (v - \Pi^\ast_{j,m_j} v) \big|}_{a, \omega_j^\circ}^2 \\
    &\le 2 |v|_{a,\omega_j^\circ}^2 + 2 {\big|\eta_j (v - \Pi^\ast_{j,m_j} v) \big|}_{a, \omega_j^\circ \cap \Omega_j^\ast}^2 \\
    &\le 2 |v|_{a,\Omega_j}^2 + 2 {\big|\eta_j (v - \Pi^\ast_{j,m_j} v) \big|}_{a, \Omega_j^\ast}^2 \le (2 + 2 / {\lambda^j_{m_j+1}}  ) |v|_{a,\Omega_j}^2\,,
  \end{align*}
  where we used that $I - H_j^\circ$ is an $a$-orthogonal projection, and the energy-minimality of the harmonic extension. Then we extend the domains in the seminorms and use Lemma~\ref{lem:stability general}. Plugging everything into~\eqref{eq:splitting} yields the result.

\end{proof}

\begin{corollary}
  The condition number of the two-level Schwarz method with the R-GenEO coarse space can be bounded as
  \begin{equation*}
    \kappa(\matr{M}^{-1}_{\mathrm{AS}, 2} \matr{A}) \le (1 + k_0) ( 2 + k_0(2 k_0 + 1) \max_{1 \le j \le N} (2 + 3 / \lambda^j_{m_j + 1} ) )\,.
  \end{equation*}
\end{corollary}

\section{Numerical results}
\begin{figure}[htpb]
  \centering
  \begin{subfigure}[b]{0.318\textwidth}
    \includegraphics[width=\textwidth]{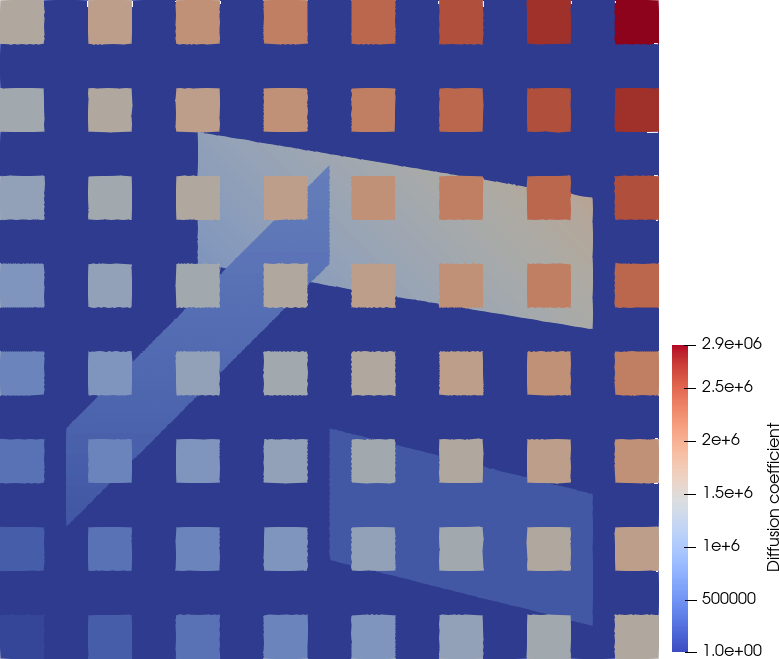}
    \caption{Diffusion coefficient}\label{fig:coeff}
  \end{subfigure}
  \hspace{2em}
  \begin{subfigure}[b]{0.27\textwidth}
    \includegraphics[width=\textwidth]{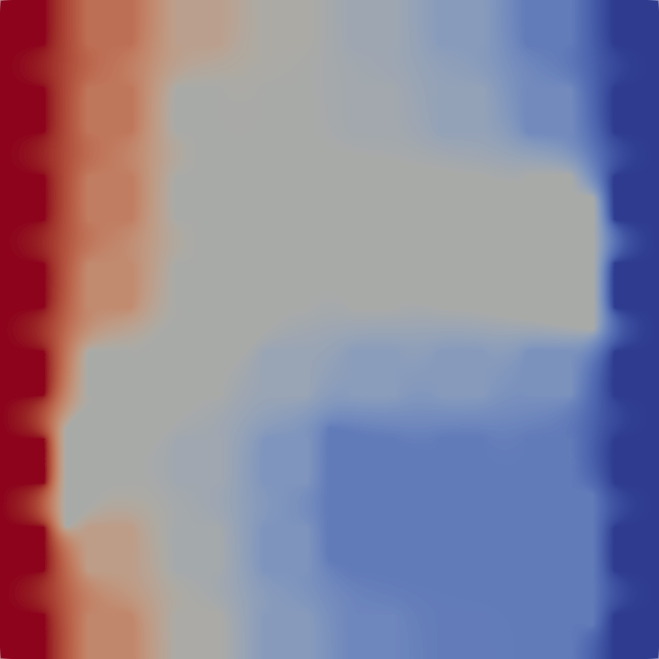}
    \caption{Finite element solution}\label{fig:solution}
  \end{subfigure}
  \caption{High-contrast coefficient and solution of the example PDE.}
\end{figure}

We consider the following variable-coefficient elliptic equation
\begin{align*}
    - \nabla \cdot (\alpha(x, y) \nabla u(x, y)) &= 0 &&\text{in } \Omega = (0,1)^2\,, \\
    u(x, y) &= 1 - x && \text{on } \{ (x,y) \in \overline{\Omega} \mid x = 0 \text{ or } x = 1 \}\,,
\end{align*}
and homogeneous Neumann boundary conditions on the rest of the boundary. The coefficient $\alpha : \Omega \rightarrow (0,\infty)$ is heterogeneous and of high-contrast ($\alpha_{\max} / \alpha_{\min} \approx 10^6$); it is visualised in Fig.~\ref{fig:coeff}. We discretise the equation using $\mathbb{Q}_1$ finite elements on a structured grid. The finite element solution is shown in Fig.~\ref{fig:solution}.

Table~\ref{table:results} reports the number of iterations to solve the resulting linear system using the conjugate gradient (CG) method preconditioned with a two-level Schwarz method using both the classical GenEO coarse space and the R-GenEO coarse space. The stopping criterion is a relative residual reduction of $10^{-10}$. We partition the domain into a regular grid of $2 \times 2$, $4 \times 4$, $8 \times 8$ or $16 \times 16$ non-overlapping subdomains and then add two layers of elements to create the overlapping subdomains. The number of elements per subdomain is $256 \times 256$ in all tests. We use a fixed number of $24$ eigenfunctions per subdomain to set up the coarse spaces.

\begin{table}[htpb]
    \setlength{\tabcolsep}{5pt}
    \centering
    \begin{tabular}{rcccccccc}\toprule
     & \multicolumn{4}{c}{GenEO} & \multicolumn{4}{c}{R-GenEO} \\ \cmidrule(lr){2-5}\cmidrule(lr){6-9}
    $N$ & \# its. & $\kappa(\matr{M}^{-1}_{\mathrm{AS}, 2} \matr{A})$ & $t_{\text{setup}}$ & $t_{\text{solve}}$ & \# its. & $\kappa(\matr{M}^{-1}_{\mathrm{AS}, 2} \matr{A})$ & $t_{\text{setup}}$ & $t_{\text{solve}}$\\\midrule
    4    & 18 & 6.49  & 4.10 & 0.37 &   19 &  6.46 & 0.96 & 0.36  \\
    16   & 27 & 13.65 & 3.70 & 0.56 &   27 & 12.71 & 1.21 & 0.52 \\
    64   & 29 & 13.40 & 4.21 & 0.70 &   30 & 13.40 & 1.26 & 0.71 \\ 
    256  & 32 & 14.58 & 4.51 & 1.14 &   31 & 14.02 & 1.49 & 1.08 \\ 
    \bottomrule
    \end{tabular}
    \caption{Comparison of the two coarse spaces. The first column: number of subdomains, `\# its.': number of CG iterations, $t_{\text{setup}}$ and $t_{\text{solve}}$: times in seconds to set up the preconditioner and to solve the linear system. The condition number $\kappa(\matr{M}^{-1}_{\mathrm{AS}, 2} \matr{A})$ was estimated after the last CG iteration using a Lanczos step.} 
    \label{table:results}
\end{table}

It can be observed that the variant based on the R-GenEO coarse space converges in a similar number of iterations while the setup-time is reduced by a factor 3.

\printbibliography
\end{document}